\documentclass[11pt]{article}

\usepackage{times}
\usepackage{graphicx}
\usepackage{amsfonts}
\usepackage{amsmath, amsthm}

\theoremstyle{plain} %% This is the default
\newtheorem{thm}{Theorem}[section]

\newtheorem{lem}[thm]{Lemma}

\theoremstyle{definition}

\theoremstyle{remark}

\newcommand{\etal}{\textit{et al.}\ }
\newenvironment{problem}[1]{\begin{quote}\textbf{#1}:}{\end{quote}}

\title{Bit retrieval: intractability and application to digital watermarking}
\author{Veit Elser \\ Department of Physics, Cornell University \\ Ithaca, NY 14853-2501 \\ USA}

\begin{document}

\maketitle

\begin{center}
\parbox{5in}{
Bit retrieval, the problem of determining a binary sequence from its cyclic autocorrelation, is a special case of the phase retrieval problem. Algorithms for phase retrieval are extensively used in several scientific disciplines, and yet, very little is known about the complexity of these algorithms or phase retrieval in general. Here we show that bit retrieval, in particular, is closely related to computations that arise in algebraic number theory and can also be formulated as an integer program. We find that general purpose algorithms from these fields, when applied to bit retrieval, are outperformed by a particular iterative phase retrieval algorithm. This algorithm still has exponential complexity and motivates us to propose a new public key signature scheme based on the intractability of bit retrieval, and image watermarking as a possible application.} 
\end{center}


\section{Introduction}
\textit{Phase retrieval} is the general problem of reconstructing a finitely sampled signal (or density in higher dimensions) from its autocorrelation. Since knowledge of the autocorrelation is equivalent to knowledge of the signal's Fourier transform modulus, phase retrieval is fundamentally underdetermined without additional information to constrain the Fourier transform phases. These constraints usually take the form of \textit{a priori} information: the signal may be known to have a particular support or distribution of values. \textit{Bit retrieval} is perhaps the simplest instance of phase retrieval, where the signal is periodic and known to take only two values. Choosing without loss of generality these values to be 0 and 1, bit retrieval seeks to find a binary sequence having a prescribed cyclic autocorrelation. For example, given the autocorrelation sequence $\alpha=[5,2,1,3,3,2,3,3,1,2]$, one solution is the binary sequence $\beta=[1,0,0,1,1,0,0,1,0,1]$.

The computational complexity of bit retrieval is largely unexplored. Zwick \etal [Zw] made the first study and were able to solve sequences up to lengths $N=64$. There is a close relationship between bit retrieval and the problem of factoring in rings of algebraic integers, specifically, the integers of the cyclotomic field of $N$th roots of unity. It is also possible to formulate bit retrieval as an integer program. While both of these subjects, algebraic number theory and integer programming, have experienced significant algorithm development in recent years, the fastest known bit retrieval algorithm still follows the principles developed in the study of phase retrieval. As described below, this algorithm has an empirically determined average-case complexity of $2^{c N}$, with $c\approx 0.22$. A point of comparison is the fact that an ordinary integer with two large factors of order $2^N$ can be factored with subexponential time complexity, specifically, $\exp{[(\log{N})^{1/3}(c\log{\log{N}})^{2/3}]}$, where $c=8/3$ [LL]. The latter problem is still considered intractable and forms the basis of public key cryptosystems [RSA]. Can the apparent intractability of bit retrieval be exploited likewise? This paper reports on a public key signature scheme as a partial response to this challenge. An application that appears to be well suited to this scheme is image watermarking.

\section{Notation and terminology}

We restrict our study to sequences of length $N$, where $N$ is an odd prime, typically greater than 200 in the applications we propose. The autocorrelation of sequences of real numbers, and more generally their \textit{convolution product}, corresponds to the standard product in the polynomial ring $\mathbb{R}[x]$. Cyclic convolutions correspond to the quotient ring $R:=\mathbb{R}[x]/\langle x^N-1\rangle$, and cyclic integer sequences form the subring $Z:=\mathbb{Z}[x]/\langle x^N-1\rangle$. 
Also of interest are the quotient rings $R/\langle\Phi_N\rangle$ and $O:=Z/\langle\Phi_N\rangle$, where $\Phi_N(x):=x^{N-1}+\cdots+x+1$ is the $N$th cyclotomic polynomial.
Since $\Phi_N(x)$ is the irreducible polynomial of $\zeta:=\exp{(\mathrm{i} 2\pi/N)}$, $O$ is isomorphic to $\mathbb{Z}[\zeta]$, the ring of integers of the cyclotomic field of $N$th roots of unity. We denote both quotient maps by the symbol $\Psi$. The \textit{rational integers} are the subring $\mathbb{Z}\subset O$.

In computations, elements of the rings $R$, $Z$ and $O$ are represented by their components with respect to a standard basis. We will use the following choice of basis elements:
\begin{eqnarray}
R, Z&:\quad& 1,x,\ldots,x^{N-1}\label{RZbasis}\\ 
O&:\quad&\quad\zeta,\ldots,\zeta^{N-1}\label{Obasis}
\end{eqnarray}
The $i$th component of an element $\alpha\in R$ is denoted $[\alpha]_i$, where $\alpha=\sum_{i=0}^{N-1}[\alpha]_i x^i$, and similarly for elements of $Z$ and $O$. 

Our bases also allow us to define the \textit{binary elements}. We say $\beta_R\in R$ is binary if $[\beta_R]_i=\pm\frac{1}{2}$ for all $0\leq i\leq N-1$ and $\Psi(\beta_R)\neq 0$. This embedding is geometrically more natural than that given in the introduction. There are exactly $2^N-2$ binary elements in $R$ and each has a distinct binary counterpart $\beta_O=\Psi(\beta_R)$ in $O$:
\begin{equation}\label{binary}
[\beta_O]_i=\left\{\begin{array}{ll}
[\beta_R]_i +\frac{1}{2} & \mbox{if $[\beta_R]_0<0$}\vspace{10pt}\cr
[\beta_R]_i -\frac{1}{2} & \mbox{if $[\beta_R]_0>0$.}\end{array}
\right.\quad 1\leq i\leq N-1
\end{equation}

The automorphisms of $O$ are given by the $N-1$ \textit{conjugate maps} defined by
$\sigma_j(\zeta):=\zeta^j$, with \mbox{$1\leq j\leq N-1$}. The collection of these maps is closely related to the
\textit{Fourier transform}. Referring the action of $\sigma_j$ on $\alpha\in O$ to the basis
(\ref{Obasis}),
\begin{equation}\label{linearSigma}
\sigma_j(\alpha)=\sum_{k=1}^{N-1}\zeta^{j\, k}[\alpha]_k\; ,
\end{equation}
we see that $\sigma:=[\sigma_1,\ldots,\sigma_{N-1}]$ can be interpreted as a linear map $O\to \mathbb{C}^{N-1}$. Since $\sum_{i=0}^{N-1}\zeta^{j\, i}=0$ for all $1\leq j\leq N-1$,  $\sigma$ is also well defined when applied to elements of $R$ and $Z$. We thus use the same notation for the set of Fourier transform components for all three rings. The statement that the automorphisms preserve multiplication in $O$, $\sigma_j(\alpha\beta)=\sigma_j(\alpha)\sigma_j(\beta)$, when written in multicomponent form as $\sigma(\alpha\beta)=\sigma(\alpha)\sigma(\beta)$, is the ``convolution theorem" of the Fourier transform. The latter, in combination with the inverse Fourier transform (see below), is the basis of an $O(N\log{N})$ multiplication algorithm (FFT) in $R$, $Z$ and $O$.

The map $\sigma_{N-1}$ corresponds to complex conjugation and will be denoted by the overbar in $O$ as well as $\mathbb{C}$: $\sigma(\overline{\alpha})=\overline{\sigma}(\alpha)$. This extends by (\ref{linearSigma}) to an action on elements of $R$ and $Z$ given by $[\overline{\alpha}]_j=[\alpha]_{N-j}$, for $1\leq j\leq N-1$, and $[\overline{\alpha}]_0=[\alpha]_0$.

The conventional Fourier transform also includes the zero-frequency component $\sigma_0\colon R\to\mathbb{R}$, where
\begin{equation}
\sigma_0(\alpha):=\sum_{i=0}^{N-1}[\alpha]_i
\end{equation}
is again to be interpreted as a linear map that extends in the obvious way to elements of $Z$.
Linear transformations $\sigma^{-1}$ and $\sigma_0^{-1}$, corresponding to the \textit{inverse Fourier transform}, are defined in the sense of the Moore-Penrose pseudoinverse:
\begin{eqnarray}
\sigma^{-1}&:=&\sigma^{\dagger}\cdot(\sigma\cdot \sigma^{\dagger})^{-1}=\frac{1}{N}\sigma^{\dagger}\label{FourierInv1}\\ 
{\sigma_0^{-1}}&:=&\sigma_0^{\dagger}\cdot(\sigma_0\cdot \sigma_0^{\dagger})^{-1}=\frac{1}{N}\sigma_0^{\dagger}\;\label{FourierInv2} ,
\end{eqnarray}
where $\cdot$ denotes matrix multiplication and $\dagger$ is the matrix adjoint. Whereas $\sigma\cdot\sigma^{-1}$ and $\sigma_0\cdot\sigma_0^{-1}$ are respectively $(N-1)\times (N-1)$ and $1\times 1$ identity matrices, the product
\begin{equation}
\pi_0:=\sigma_0^{-1}\cdot\sigma_0=\frac{1}{N}\left({\begin{array}{ccc}
1&\cdots&1\\
\vdots& &\vdots\\
1&\cdots&1
\end{array}}\right)
\end{equation}
is the projector to the ideal $\langle\Phi_N\rangle$ in $R$. Similarly, $\sigma^{-1}\cdot\sigma$ is the projector onto the orthogonal complement, $R_\perp$, where orthogonality is with respect to the \textit{Euclidean norm}:
\begin{eqnarray}
\|\alpha\| &:=& \left(\sigma_0(\alpha)^2+\overline{\sigma}(\alpha)\cdot\sigma(\alpha)\right)/N\label{quadForm1}\\
&=& \alpha^{\mathrm{t}}\cdot\pi_0\cdot\alpha+\alpha^{\mathrm{t}}\cdot(1-\pi_0)\cdot\alpha\label{quadForm2}\\
&=& \alpha^{\mathrm{t}}\cdot\alpha\; .\label{quadForm3}
\end{eqnarray}
The Euclidean norm for elements $\alpha\in R_\perp$,
\begin{equation}\label{Onorm}
\|\alpha\|_{\perp}:=\overline{\sigma}(\alpha)\cdot\sigma(\alpha)/N\; ,
\end{equation}
is also the appropriate norm in the quotients $R/\langle\Phi_N\rangle$ and $O$. 
Some of the interest in studying binary elements derives from the fact that all binary $\beta\in R$ have the same Euclidean norm, $\|\beta\|=N/4$.

The \textit{autocorrelation} $\alpha$ of an element $\beta\in R, Z, O$ is given by $\alpha=\beta\overline{\beta}$ and has real, nonnegative Fourier transform components:  $\sigma(\beta\overline{\beta})=\sigma(\beta)\overline{\sigma}(\beta)$, and also $\sigma_0(\beta\,\overline{\beta})=\sigma_0(\beta)^2$ for $\beta\in R, Z$. The autocorrelation of an element $\beta$ is therefore equivalent to the information in its Fourier transform modulus, and recovering $\beta$ from its autocorrelation corresponds to ``retrieving its phases". Autocorrelations, and more generally, elements with the property $\overline{\alpha}=\alpha$, form the real subrings $\hat{R}$, $\hat{Z}$ and $\hat{O}$. If $\beta_O\in O$ is binary and $\beta_R$ is its binary counterpart in $R$, then the corresponding autocorrelations $\alpha_O=\beta_O\overline{\beta}_O$ and $\alpha_R=\beta_R\overline{\beta}_R$ are related by
\begin{equation}
\begin{split}
&{[\alpha_R]_0}=\frac{N}{4}\\
&{[\alpha_R]_i}=[\alpha_O]_i+\frac{N}{4}\; ,\quad 1\leq i\leq N-1\; .\label{autoRandO}
\end{split}
\end{equation}
 
A binary element $\beta\in O$ is said to be \textit{perfect} if its autocorrelation is a rational integer, that is, $\beta\overline{\beta}\in \mathbb{Z}$. The Fourier transform components of a perfect $\beta$ have constant modulus, since $\beta\overline{\beta}=\sigma_j(\beta\overline{\beta})=|[\sigma(\beta)]_j|^2$. For any $N$, $\beta=1$ is perfect; a less trivial example, for $N=3$, is the binary element $\beta=1+\zeta$.

The norm $\mathcal{N}(\alpha)\in\mathbb{Z}$ of an element $\alpha\in O$ is defined by
\begin{equation}\label{norm}
\mathcal{N}(\alpha):=\prod_{j=1}^{N-1}\sigma_j(\alpha)=\prod_{j=1}^{(N-1)/2}|\sigma(\alpha)_j|^2\;,
\end{equation}
and has the interpretation of the index in $O$ of the principal ideal $\alpha\, O$.

\section{Bit retrieval}

A generalization of the problem posed in the introduction is the following:
\begin{problem}{B$_1$}
Given $\alpha\in O$ and the knowledge that $\alpha=\beta_1\, \beta_2$ where $\beta_1$ and $\beta_2$ are binary, find a particular such pair $\beta_1$ and $\beta_2$.
\end{problem}
The security of the proposed signature scheme relies on the intractability of two related problems:
\begin{problem}{B$_2$}
Given $\alpha\in \hat{O}$ and the knowledge $\alpha=\beta\overline{\beta}$ where $\beta$ is binary, find such a $\beta$.
\end{problem}
\begin{problem}{B$_3$}
Given a finite set $A\subset O$ and the knowledge that some binary element $\beta\in O$ divides every $\alpha\in A$, find such a $\beta$.
\end{problem}

In the ring of rational integers these problems correspond to factorization (\textbf{B$_1$}), finding the square root of a perfect square (\textbf{B$_2$}), and obtaining the GCD of a set of integers (\textbf{B$_3$}). Of these, only factorization remains intractable, the square root and GCD being computed efficiently by Newton's and Euclid's algorithms respectively. The failure of unique factorization in $O$, already for $N\geq 23$ [MM], implies that a Euclidean algorithm is not available for efficiently solving \textbf{B$_3$} in these rings. Although \textbf{B$_2$} and \textbf{B$_3$} are clearly easier than \textbf{B$_1$}, what makes the ring $O$ attractive is that even the former problems appear to be intractable when the size of the problem corresponds to $N$, rather than the sizes of the rational integers in the specification (coefficients in the standard basis). It is also for this reason that we restrict the unknown factors or ``square roots" to be binary.

Clearly problem \textbf{B$_2$} becomes easy when the density of 0's in the binary element $\beta$ is either very large or very small. Rankenburg [R] shows that the symmetric case $\beta=\overline{\beta}$ also represents an easy instance of \textbf{B$_2$}. For symmetric $\beta$ the unknown phases are either $0$ or $\pi$, and in particular, one of the following equations holds: $\sigma_1(\beta)=\pm\sqrt{\sigma_1(\alpha)}$. Solving either equation for the set of unknown binary components of $\beta$ is equivalent to solving subset-sum problems of arbitrarily low density (see section \ref{algNumAlgorithm}), and methods based on lattice basis reduction [LO] provide a polynomial-time algorithm.  

While the algebraic statements of the bit retrieval problems above seem natural, the most efficient known algorithm for solving \textbf{B$_2$}, in particular, is entirely non-algebraic. For this algorithm (section \ref{diffMapAlgorithm}), as well as integer programming methods (section \ref{intProgAlgorithm}), what matters is the following formulation as a geometric \textit{feasibility problem} in the ring $R$. We recall that autocorrelations of corresponding binary elements in the rings $O$ and $R$ are simply related by (\ref{autoRandO}).

Consider two subsets of $R$: the hypercube
\begin{equation}\label{hypercube}
B:=\left\{\beta\in R\colon [\beta]_i=\pm\frac{1}{2},\; 0\leq i\leq N-1\right\}\; ,
\end{equation}
and for any $\alpha\in \hat{R}$, the set
\begin{equation}\label{torus}
T_\alpha:=\{\beta\in R\colon \beta\overline{\beta}=\alpha \}\; .
\end{equation}
The restatement of \textbf{B$_2$} as a feasibility problem is then:
\begin{problem}{B$^\prime_2$}
Given $\alpha\in \hat{R}$, known to be the autocorrelation of a binary element, find $B\cap T_\alpha$.
\end{problem}
When characterized by its Fourier transform, the set $T_\alpha$ is recognized as a pair of $(N-1)/2$ dimensional tori. Let $\beta\in T_\alpha$, then the definition (\ref{torus}) implies
\begin{equation}\label{torusRadii}
\begin{split}
&\sigma_0(\beta)=\pm\sqrt{\sigma_0(\alpha)}\\
&|\sigma_j(\beta)|=\sqrt{\sigma_j(\alpha)}\; ,\quad 1\leq j\leq (N-1)/2\; ,
\end{split}
\end{equation}
with no further constraints required on the remaining components because of complex conjugation symmetry. Using the linearity of the Fourier transform it is straightforward to show that the convex hull of $T_\alpha$ is given by
\begin{equation}\label{torusHull}
h(T_\alpha):=\left\{\beta\in R\colon |\sigma_i(\beta)|\leq\sqrt{\sigma_i(\alpha)},\; 0\leq i\leq (N-1)/2\right\}\; .
\end{equation}
Since convex relaxations of constraints typically simplifies feasibility problems, we also consider the convex hull of the hypercube,
\begin{equation}
h(B):=\left\{\beta\in R\colon |[\beta]_i|\leq\frac{1}{2},\; 0\leq i\leq N-1\right\}\; .
\end{equation}
The convex relaxations that apply to problem \textbf{B$^\prime_2$} are summarized in the following:
\begin{thm}\label{feasibilityFormulations}
Let $\alpha\in \hat{R}$ be the autocorrelation of a binary element; then
\begin{equation}
B\cap T_\alpha=h(B)\cap T_\alpha=B\cap h(T_\alpha)
\end{equation}
\end{thm}
\begin{proof}
The equality of these sets follows from the observation that if $\beta\in h(B)$ then $\|\beta\|\leq N/4$ and equality requires $\beta\in B$. Similary, if $\tau\in h(T_\alpha)$, then 
\begin{equation}
\sum_{i=0}^{N-1}|\sigma_i(\tau)|^2\leq\sum_{i=0}^{N-1}\sigma_i(\alpha)=\sum_{i=0}^{N-1}|\sigma_i(\beta)|^2\; ,
\end{equation}
where $\beta$ is some binary element. Thus
$\|\tau\|\leq \|\beta\|=N/4$ and equality implies $\tau\in T_\alpha$. Now suppose $\gamma\in h(B)\cap T_\alpha$; then since $\gamma\in T_\alpha$ we know $\|\gamma\|=N/4$. On the other hand, since $\gamma\in h(B)$, this norm is possible only if in fact $\gamma\in B$. The same argument shows that $B\cap h(T_\alpha)=B\cap T_\alpha$.
\end{proof}

\subsection{Uniqueness in bit retrieval}

For the digital signature scheme considered in section \ref{signature}, which derives its security from the conjectured intractability of bit retrieval, there is no requirement that the solutions to any of problems \textbf{B$_1$}, \textbf{B$_2$}, or \textbf{B$_3$} be unique. As described in more detail in section \ref{security}, this scheme only requires, more generally, that it is difficult  to find any solution with small Euclidean norm. It is interesting nevertheless, to ask what varieties of non-uniqueness can occur in bit retrieval. Our remarks here will address problem \textbf{B$_2$}.

Clearly if $\beta\in O$ solves \textbf{B$_2$} then so does $\overline{\beta}$. This together with the statement expressed in the following lemma characterizes the symmetries inherent in bit retrieval.
\begin{lem}
If $\beta\in O$ and $\beta\gamma\in O$ are two solutions of an instance of \textbf{B$_2$}, then $\gamma=\pm\zeta^k$ for some $k$.
\end{lem}
\begin{proof}
Since both solutions must have the same autocorrelation, $\gamma\overline{\gamma}=1$. This implies $(\mathcal{N}(\gamma))^2=1$ and we infer that $\gamma$ is a unit. Kummer's lemma can now be used to rewrite the autocorrelation of $\gamma$ with the result $\gamma^2=\zeta^k$ for some $k$. This shows that $\gamma$ is a $2N$-th root of unity, as asserted. 
\end{proof}

Beyond the symmetries that apply to any solution, problem instances can suffer from special forms of non-uniqueness. One of these has a counterpart in crystallographic phase retrieval [PS] and applies when a solution is a product $\beta=\beta_1\beta_2$, and neither factor is of the form $\pm\zeta^k$. It may then happen that $\beta^\prime=\beta_1\overline{\beta_2}$ is also binary and not related to $\beta$ by one of the symmetries discussed above. Since $\beta^\prime$ has the same autocorrelation as $\beta$, it also solves \textbf{B$_2$}. An example of this mechanism for $N=13$ arises for the autocorrelation $\beta\overline{\beta}=3$. From the factors $\beta_1=1+\zeta^2+\zeta^7$ and $\beta_2=1+\zeta^3+\zeta^4$ one obtains $\beta=-\zeta-\zeta^8-\zeta^9-\zeta^{12}$ and $\beta^\prime=-\zeta-\zeta^5-\zeta^6-\zeta^8$ as the two binary solutions. Instances with non-unique solutions, such as this one, become very rare as $N$ increases. In a set of experiments with $23\leq N\leq 53$, random binary $\beta$ were drawn from the uniform distribution using a pseudo-random number generator. When the autocorrelation of $\beta$ was given to the difference map algorithm (see below), the solution $\beta^\prime$ was compared with $\beta$. The fraction of solutions $\beta^\prime$ not symmetry-related to $\beta$ was found to decrease rapidly with $N$, as shown in Table 1.

\begin{table}
\begin{center}
\begin{tabular}{c|c|c|c|c|c|c|c}
  $N=23$ & 29 & 31 & 37 & 41 & 43 & 47 & 53\\
 $0.044$	& $0.024$ & $0.019$ & $6.1\times 10^{-3}$ & $2.4\times 10^{-3}$ & $1.4\times 10^{-3}$ & $4.7\times 10^{-4}$ & $1.5\times 10^{-4}$\\

\end{tabular}
\caption{Probability of non-uniqueness in bit retrieval}
\end{center}
\end{table}

\subsection{Facts concerning the norm}

With $N$ fixed, what characterizes hard instances of bit retrieval? The norm $\mathcal{N}(\beta)$ of the secret binary element $\beta\in O$ is a natural candidate and in fact establishes a connection with the subject of cyclic difference sets.
\begin{thm}\label{normBoundThm}
Let $\beta\in O$ be binary, then $\mathcal{N}(\beta)\leq (\frac{N+1}{4})^{\frac{N-1}{2}}$ and equality holds only if $\beta$ is perfect.
\end{thm}
\begin{proof}
Let $\beta_R$ be the binary counterpart in $R$ of a binary element $\beta\in O$ (see (\ref{binary})); then $[\beta_R]_i=\pm\frac{1}{2}$ for $0\leq i\leq N-1$. Let $\mu_j:=|[\sigma(\beta)]_j|^2=|[\sigma(\beta_R)]_j|^2$ denote the squares of the corresponding Fourier moduli. Using expressions (\ref{quadForm1}, \ref{quadForm3}) for the Euclidean norm and (\ref{norm}) for the algebraic norm, we have:
\begin{eqnarray}
\sigma_0(\beta_R)^2+\sum_{j=1}^{N-1}\mu_j &=& N\,\beta_R\cdot \beta_R = \frac{N^2}{4}\; ,\\
\prod_{j=1}^{N-1}\mu_j &=& \mathcal{N}(\beta)^2\; .
\end{eqnarray}
Applying the arithmetic-geometric mean inequality to the numbers $\mu_j$ we obtain:
\begin{equation}
\mathcal{N}(\beta)\leq \left[\frac{1}{N-1}\left(\frac{N^2}{4}-\sigma_0(\beta_R)^2\right)\right]^{\frac{N-1}{2}}\; .
\end{equation}
Since $\beta_R$ has an odd number of $\pm\frac{1}{2}$ components, $\sigma_0(\beta_R)^2\geq\frac{1}{4}$ and the stated bound on $\mathcal{N}(\beta)$ follows. Equality of the arithmetic and geometric means requires that the (squared) Fourier moduli $\mu_j$ are equal, and this is one way of characterizing a perfect $\beta$. 
\end{proof}

We will refer to binary elements that achieve the upper bound in theorem \ref{normBoundThm} as \textit{Hadamard} because of their direct relationship to Hadamard cyclic difference sets. More generally [Ba], a \textit{cyclic difference set} can be defined in terms of the cyclic group $G$ of order $N$ acting on binary elements $\beta\in R$ with generator $g: \beta\mapsto x \beta$. Defining a subset of $G$ by $D:=\{g^i\colon [\beta]_i>0, 0\leq i \leq N-1\}$, we can ask if it is possible for every nonidentity element of $G$ to appear exactly $\lambda$ times in the set $\{d_1{d_2}^{-1}\colon d_1, d_2\in D\}$. If this is the case, and the cardinality of $D$ is $k$, then 
\begin{equation}\label{diffSet1}
(N-1)\lambda+k=k^2\; ,
\end{equation}
and $D$ is declared a cyclic difference set with parameters $(N, k, \lambda)$. The binary $\beta$ which defines such a difference set will then satsify $(\beta+\frac{1}{2}\Phi_N)(\overline{\beta}+\frac{1}{2}\Phi_N)=[k, \lambda, \lambda, \ldots, \lambda]$, that is, $\beta$ will be perfect since $\Psi(\beta\overline{\beta})= k-\lambda$ (a rational integer). A \textit{Hadamard} cyclic difference set maximizes $k-\lambda$ to the maximum value consistent with the norm bound from theorem \ref{normBoundThm}: 
\begin{equation}\label{diffSet2}
k-\lambda = \frac{N+1}{4}\; .
\end{equation}
From (\ref{diffSet1}) and (\ref{diffSet2}) one obtains the Hadamard cyclic difference set parameters $(N, \frac{N-1}{2}, \frac{N-3}{4})$, which evidently require that $N\equiv 3\bmod{4}$.

There is a simple construction of Hadamard cyclic difference sets for any prime $N$ of the form $4m+3$ [Ba]; the formula for the corresponding binary $\beta\in O$ is given by:
\begin{equation}
[\beta]_i=\frac{1-(i|N)}{2}\quad (1\leq i \leq N-1)\; ,
\end{equation}
where the Legendre symbol $(i|N)$ equals $1$ whenever $i$ is a square in the finite field of order $N$, and $-1$ otherwise. For certain special values of $N$, such as $N=2^m-1$ and $N=4m^2+27$, other constructions of Hadamard cyclic difference sets are known [Ba]. An example of a Hadamard integer for $N=7$ is $\beta=1+\zeta^2+\zeta^3$.

The norms of ``random" binary integers are typically significantly smaller than the norm of a Hadamard integer. This is made precise in the following theorem.
\begin{thm}\label{normRand}
Let $\beta\in O$ be treated as a discrete random variable with uniform distribution on the set of binary integers; then as $N\to \infty$ the random variable $S:=\log{\mathcal{N}(\beta)}$ has expectation value
\begin{equation}
\mathrm{E}(S)=\frac{1}{2}\left(\log{\left(N/4\right)}-\gamma\right)N\label{E(S)}\; ,
\end{equation}
where $\gamma=0.577215\ldots$ is Euler's constant.
\end{thm}
\begin{proof}
Define the random variables $z_j:=\sigma_j(\beta)\in\mathbb{C}$, $1\leq j\leq \frac{N-1}{2}$. Each $z_j$ is the sum of $N$ independent two-valued random variables, for which the Lindeberg criterion [Bi] is easily verified. Thus as $N\to \infty$ each $z_j=x_j+\mathrm{i}y_j$ is normally distributed in $\mathbb{C}$ with distribution
\begin{equation}
\mathrm{P}(z_j)dx_j\, dy_j=\frac{4}{\pi N}\exp{\left(-4|z_j|^2/N\right)}\, dx_j\, dy_j\; .
\end{equation}
The desired expectation value may now be calculated as follows:
\begin{eqnarray}
\mathrm{E}(S)&=&\mathrm{E}\left({\textstyle \sum_{j=1}^{(N-1)/2}}\log{|z_j|^2}\right)\\
&\sim&\frac{N}{2}\int \log{|z|^2}\,\mathrm{P}(z)\,dx\, dy\quad (N\to\infty)\\
&=&\frac{N}{2}\int_0^\infty \log{(t N/4)}\,e^{-t}\, dt\; ,
\end{eqnarray}
and the stated result (\ref{E(S)}) follows.
\end{proof}

The norm of a random binary integer is thus smaller by a factor of order $\exp{(-\gamma N/2)}$, relative to the norm of a Hadamard integer. Below it is speculated that this may account for the fact that the difference map algorithm typically requires many more iterations for the retrieval of a Hadamard instance. Although the difference map algorithm is non-algebraic and works with the geometric formulation \textbf{B$^\prime_2$}, the norm is still relevant because of the fact expressed by the following theorem. 

\begin{thm}
Let $\beta_R$ be an embedding of $\beta\in O$ in $R$,
\begin{eqnarray}
{[\beta_R]_0} &=&r\\
{[\beta_R]_j} &=&[\beta]_j+r\; ,\quad 1\leq j\leq N-1\; ,
\end{eqnarray}
where $r\in \mathbb{R}$ is arbitrary. Then
\begin{equation}
\mathrm{vol}(T_\alpha)=2\left(\frac{8\pi^2}{N}\right)^{\frac{N-1}{4}}\sqrt{\mathcal{N}(\beta)}\; ,
\end{equation}
where $T_\alpha$ is the torus defined in (\ref{torus}) and specified by $\alpha=\beta_R\overline{\beta}_R$.
\end{thm}

\begin{proof}
Consider a point $\tau\in T_\alpha$. From (\ref{torusRadii}) we infer
\begin{eqnarray}
\sigma_0(\tau) &=& \pm |\sigma_0(\beta_R)|\\
\sigma_j(\tau) &=& |\sigma_j(\beta_R)|\exp{\mathrm{i}\phi_j} = |\sigma_j(\beta)|\exp{\mathrm{i}\phi_j}\; ,
\end{eqnarray}
where the angles $\phi_j$ for $j=1,\ldots, \frac{N-1}{2}$ are arbitrary and related to the others by $\phi_j=-\phi_{N-j}$. This shows that topologically $T_\alpha$ comprises two smooth tori of dimension $\frac{N-1}{2}$. The angles $\phi_j$ serve as convenient coordinates in the explicit representation for a general point $\tau\in T_\alpha$:
\begin{equation}
\tau = \pm \sigma_0^{-1} \cdot |\sigma_0(\beta_R)|+\sigma^{-1}\cdot |\sigma(\beta)|\exp{\mathrm{i}\phi}\; .
\end{equation}
The computation of the volume is now an elementary exercise in calculus and leads directly to the quoted value.
\end{proof}

\section{Algorithms}

Bit retrieval falls within the scope of at least three algorithmic frameworks: (i) algebraic number theory, (ii) integer programming, and (iii) phase retrieval. We describe below all three as they apply to problem \textbf{B$_2$}, or its geometrical formulation \textbf{B$^\prime_2$}. Problems \textbf{B$_1$} and \textbf{B$_3$} are almost indistinguishable from \textbf{B$_2$} within the algebraic approach, whereas the integer programming and phase retrieval techniques first require geometrical reformulations of \textbf{B$_1$} and \textbf{B$_3$} before these methods can be applied to them.

\subsection{Algebraic number theory}\label{algNumAlgorithm}

In the algebraic approach the secret binary integer ($\beta$ in problems \textbf{B$_2$} and \textbf{B$_3$}, $\beta_1$ or $\beta_2$ in \textbf{B$_1$}) is first identified by the principal ideal it generates in $O$: $I=\langle\beta\rangle$. This task is relatively easy and almost insignificant in comparison to the subsequent task of finding a binary generator of $I$. There are algorithms [Co] that take as input the generators of an ideal $I$ and return a single generator $\gamma$ if $I$ is found to be principal. This would appear to be a good technique, since the desired binary generator $\beta$ can then be expressed in the form $\beta= u\gamma$, where $u$ is a unit. However, algorithms for principal ideal testing require information about the class group of $O$, making this approach prohibitive already for $N\ge 67$ [Bu]. An alternative, used in the algorithm below, is to work only with the lattice structure of $I$ and seek a binary element $\beta^\prime\in I$ without the guarantee that $\langle\beta^\prime\rangle=I$. Since there are so few binary elements in $I$, a practical approach is to enumerate them completely using the Fincke-Pohst algorithm [FP] and thereby discover the particular element that generates $I$. The complexity of the algebraic approach is thus determined by the complexity of an associated lattice search problem.

An example with $N=23$ should serve as a substitute for a formal specification of the algorithm. The identity of the secret binary $\beta$ is contained in its autocorrelation $\alpha=\beta\overline{\beta}$, say
\begin{equation}
\alpha=-[5, 7, 4, 5, 7, 7, 5, 6, 8, 6, 6, 6, 6, 8, 6, 5, 7, 7, 5, 4, 7, 5]\; ,
\end{equation}
or in products $\gamma_1=\beta\beta_1$, $\gamma_2=\beta\beta_2$, etc. Suppose we are given just two:
\begin{eqnarray}
\gamma_1&=&[3, 0, 0, 2, 0, -1, -1, 1, 0, -2, 0, 0, 1, 2, 0, 3, 2, 0, 2, 2, 2, -1]\; ,\\
\gamma_2&=&[0, 2, 0, -1, 0, 1, 0, 1, 0, 0, -1, -1, 0, 0, -1, 1, 1, 1, 1, 0, 0, 0]\; .
\end{eqnarray}
Using efficient algorithms (see [Co]) the ideals generated by $\alpha$, $\gamma_1$ and $\gamma_2$ can be factored into prime ideal factors with the following result:
\begin{equation}\label{primeFactors}
\langle\alpha\rangle = I_1 I_2 I_3 I_4\quad\langle\gamma_1\rangle = I_1 I_4 I_5\quad\langle\gamma_2\rangle = I_1 I_4 I_6 I_7
\end{equation}
\begin{equation}
I_1 = \overline{I}_2 = \langle 47, 15+\zeta\rangle \quad I_3 = \overline{I}_4 = \langle 5843, 1833+\zeta\rangle
\end{equation}
\begin{equation}
I_5 = \langle 174157, 61966+\zeta\rangle\quad
I_6 = \langle 47, 13+\zeta\rangle\quad
I_7 = \langle 1979, 152+\zeta\rangle
\end{equation}
As an example of an instance of problem \textbf{B$_1$} we would be given only $\gamma_1$, say, and the factorization (\ref{primeFactors}) would provide us with eight candidate factorizations of $\langle\beta\rangle$. This includes the rare possibility that $\beta$ is a unit. The number of trial factorizations to explore will almost always be small, and this is especially the case for the other two bit retrieval problems. In problem \textbf{B$_3$} we have factorizations for both $\langle\gamma_1\rangle$ and $\langle\gamma_2\rangle$, giving only four possible factorizations of $\langle\beta\rangle$. Moreover, the random origins of $\beta_1$ and $\beta_2$, say in a digital signature scheme, would imply $\langle\beta\rangle=I_1 I_4$ with high probability. Finally, in problem \textbf{B$_2$} we know that $\alpha$ decomposes into a complex conjugate pair giving only two possibilities to consider, $\langle\beta\rangle=I_1 I_3$ and $\langle\beta\rangle=I_1 I_4$.

For each candidate factorization, the number of which will be small, another standard algorithm [Co] returns the lattice basis of the corresponding ideal product in Hermite normal form. Given this basis we can in principle determine if the lattice contains a nonzero binary vector. From experiments with ideals generated by random binary elements we find that with high probability the Hermite normal form basis has the following simple form:
\begin{equation}\label{lattice}
v_j := a_j\zeta+\zeta^j\quad (1\leq j\leq N-1)\; .
\end{equation}
This is also the case for the correct factorization choice in our example, $\langle\beta\rangle=I_1 I_4$, where
\begin{equation}
\begin{split}
a = &[274620, 218518, 159293, 98597, 171309, 37690, 214991, 11132, 50442, 252742, 78333,\\
 & 231057, 55808, 42203, 207268, 79601, 242822, 193340, 248383, 212667, 72735, 58266]\; .
\end{split}
\end{equation}
We note that $a_1+1=274621=\mathcal{N}(\beta)$. In general, lattices of high index are unlikely to contain any nonzero binary vectors, in particular, the secret $\beta$. Given a ``random" lattice of index $\mathcal{N}(\beta)$ one expects to find $(2^{N-1}-1)/\mathcal{N}(\beta)$ binary vectors, a number which vanishes with $N$ as $(N e^{-\gamma}/16)^{-N/2}$ using the asymptotic result of theorem \ref{normRand}. We may therefore conclude that an exhaustive search for nonzero binary vectors in the lattice generated by the $v_j$ will either yield no results, as in fact happens when the wrong factorization $\langle\beta\rangle=I_1 I_3$ is tried, or will produce just the desired solutions $\beta\zeta^i$, $1\leq i\leq N$. Any binary element $\beta^\prime$ produced by the search must be tested against the given autocorrelation $\alpha$ since, as an element $\beta^\prime\in\langle\beta\rangle$, we only have the guarantee that $\alpha$ divides $\beta^\prime\overline{\beta^\prime}$. This does not pose a problem in practice since $\beta^\prime\overline{\beta^\prime}\neq\alpha$ implies $\mathcal{N}(\beta^\prime)\ge\mathcal{N}(\beta)$, corresponding to an even smaller expected number of binary vectors with the incorrect autocorrelation. For the example above, in fact, the search found only the true solution
\begin{equation}\label{pi23}
\beta=[1, 1, 0, 0, 1, 0, 0, 1, 0, 0, 0, 0, 1, 1, 1, 1, 1, 1, 0, 1, 1, 0]
\end{equation}
and its 22 multiples with powers of $\zeta$.  

For lattice bases of the form (\ref{lattice}), the problem of finding a binary vector is closely related to a \textit{subset sum} problem. Let $A=\{a_2, a_3,\dots a_{N-1}\}$, then finding a binary vector is equivalent to finding a subset $A^\prime\subset A$ with sum $\Sigma(A^\prime)$, such that $\Sigma(A^\prime)\in \{0,1\}\pmod{\mathcal{N}(\beta)}$. Because the subset sum problem is known to be NP-complete [GJ], this approach to bit retrieval cannot guarantee a polynomial-time solution. However, by expressing the subset sum problem as a shortest lattice vector problem, Lagarias and Odlyzko [LO] showed that instances with sufficiently small density $d$ can be solved efficiently using lattice basis reduction algorithms, where
\begin{equation}
d:=\frac{|A|}{\max_{a\in A}(\log_2{a})}\; .
\end{equation}
Evaluating this for bit retrieval instances, where $|A|=N-2$ and $a<\mathcal{N}(\beta)$ for all $a\in A$, we obtain
\begin{equation}\label{densBound}
d>\frac{N-2}{\log_2{\mathcal{N}(\beta)}}\stackrel{N\to\infty}{\sim}\frac{\log{4}}{\log{(N/4)}-\gamma}\; ,
\end{equation}
using the result of theorem \ref{normRand}.
The bound (\ref{densBound}) violates the criterion found by Lagarias and Odlyzko, who showed that $d$ must be no greater than $O(1/N)$ in order for the LLL polynomial-time basis reduction algorithm [LLL] to succeed in solving the subset sum problem.

The Fincke-Pohst nearest vector algorithm [FP] would appear to be the best technique for finding a binary vector since it guarantees a solution regardless of density while also taking advantage of LLL basis reductions. When given the generators $v_j$ and target vector $[\frac{1}{2},\cdots,\frac{1}{2}]$, this algorithm returns all binary vectors in the lattice generated by the $v_j$. Table 2 gives running times for the \texttt{kant} [K] implementation of this algorithm on bit retrieval instances up to $N=41$. All instances were generated by taking the leading $N-1$ base 2 digits of $\pi=11.001\ldots$ as the components of the secret binary integer $\beta\in O$ in the standard basis. These same ``$\pi$-sequence" instances, $\beta=\pi_N$, were used to test the other two algorithms discussed below. The solution given in (\ref{pi23}) is $\pi_{23}$.

\begin{table}
\begin{center}
\begin{tabular}{c|r|r|r|r|r|r|r|r}

algorithm  & $N=23$ & 29 & 31 & 37 & 41 & 43 & 47 & 53 \\
\hline

algebraic number theory (\texttt{kant4}) & 0.8 (sec) & 9.9 & 31 & 3800 & 62000 & * & * & * \\
integer programming (\texttt{bonsaiG}) & 0.2 (sec) & 27 & 7.2 & 79 & 8000 & 4300 & 11000 & * \\
phase retrieval (difference map) & $<0.1$ (sec) & $<0.1$ & $<0.1$ & $<0.1$ & 0.4 & 1.1 & 0.5 & 2.9 \\
\hline
\end{tabular}
\parbox{5in}{\caption{Timing results for three bit retrieval algorithms on $\pi$-sequence instances for software running on a single 1.67 GHz Athlon processor (* time limit exceeded).}}

\end{center}
\end{table}

It is probably no coincidence that the long running times for $N>31$ coincide with the relatively abrupt onset of the LLL algorithm's inability to discover generators for ideals $\langle\beta\rangle$ when given a lattice basis in Hermite normal form. Results for the latter problem are shown in Table 3. In these experiments LLL reduction was applied to the Hermite normal form basis of the principal ideal generated by a random binary element $\beta\in O$. A successful instance of principal ideal discovery was declared if one of the reduced basis elements $v^{\prime}_j$ satisfied $\mathcal{N}(v^{\prime}_j)=\mathcal{N}(\beta)$. From the results in Table 3 we see that the success rate vanishes rapidly with increasing $N$, beginning at about $N=31$.

\subsection{Integer programming}\label{intProgAlgorithm}
The form of the feasibility problem \textbf{B$^\prime_2$} that is most amenable to the techniques of integer programming is that given in theorem \ref{feasibilityFormulations}, of finding an element in the intersection $B\cap h(T_\alpha)$. Although $h(T_\alpha)$ is convex, standard integer programming algorithms based on linear relaxations also require that this set be defined by linear inequalities. We therefore make the further relaxation of replacing $h(T_\alpha)$, geometrically a product of disks, by a product of squares (and one interval):
\begin{equation}\label{torusRelax}
\begin{split}
sh(T_\alpha):=&\left\{\beta\in R\colon |\sigma_0(\beta)|\leq \sqrt{\sigma_0(\alpha)}\; , \right. \\
 &\left. |\Re{(\sigma_j(\beta))}|\leq \sqrt{\sigma_j(\alpha)}\; ,\; |\Im{(\sigma_j(\beta))}|\leq \sqrt{\sigma_j(\alpha)}\; ,\; 1\leq j\leq N-1 \right\}\; .
\end{split}
\end{equation}
Since $h(T_\alpha)\subset sh(T_\alpha)$, all bit retrieval solutions are contained in $B\cap sh(T_\alpha)$. Although we cannot rule out the possibility $B\cap sh(T_\alpha)\neq B\cap h(T_\alpha)$, this is a concern only if the relaxed problem admits too many additional solutions. Experiments show that in fact this is not the case: only bit retrieval solutions were found in all the instances studied.

In standard linear programming notation, the feasibility problem for $B\cap sh(T_\alpha)$ is expressed as:
\begin{problem}{find}
$b \in \{-\frac{1}{2},\frac{1}{2}\}^{N}$
\end{problem}
\begin{problem}{such that}
$|C\cdot b|\leq a\quad$ and $\quad |S\cdot b|\leq a$
\end{problem}
\begin{problem}{where}
$a_i = \sqrt{\sigma_i(\alpha)}\quad C_{ij} = \cos{(2\pi i j/N)}\quad S_{ij} = \sin{(2\pi i j/N)}\quad(0\leq i, j\leq N-1)$
\end{problem}
This linear program comprises exactly $2N$ independent and nontrivial inequalities for $N$ binary variables. Somewhat unusual is the fact that the coefficient matrices have nearly unit density. Solution times for the general-purpose solver \texttt{bonsaiG} [Ha] on the $\pi$-sequence instances are given in Table 2. Over the limited range studied, it appears the performance of the integer programming algorithm is somewhat better than that of the algebraic number theory based algorithm.

\begin{table}[t]

\begin{center}
\begin{tabular}{c|c|c|c|c|c|c|c}

 $N=29$ & 31 & 37 & 41 & 43 & 47 & 53 & 59\\

$0.923$	& $0.851$ & $0.504$ & $0.232$ & $0.158$ & $0.070$ & $0.011$ & $0.002$\\

\end{tabular}
\caption{Success rate of principal ideal discovery by LLL basis reduction}
\end{center}
\end{table}

\subsection{Phase retrieval}\label{diffMapAlgorithm}

Because the constraints in phase retrieval are typically nonconvex, very different solution strategies have evolved to solve these problems. Although not true algorithms in a strict sense, with a bounded running time, these methods are very successful and are not likely to be replaced by more rigorously defined algorithms in the near future. Here we apply a general purpose phase retrieval method, the \textit{difference map} [E1], to problem \textbf{B$^\prime_2$}. The difference map applies to the general feasibility problem of finding an element in $A\cap B$, where $A$ and $B$ are arbitrary sets in a Euclidean space. Practical implementations of the difference map are limited to situations where the projectors $\Pi_A$ and $\Pi_B$, to respectively the sets $A$ and $B$, can be computed efficiently. A brief description of the method is given in the Appendix.

We choose for our two sets the torus $T_\alpha$ and hypercube $B$ (as instances of the general sets $A$ and $B$ of the Appendix); experimentation indicates there is no advantage in using either of the convex relaxations given in theorem \ref{feasibilityFormulations}. The projectors $\Pi_{T_\alpha}$ and $\Pi_B$ are maps $R\to R$ where
\begin{equation}
\Pi_{T_\alpha} := \sigma_0^{-1}\cdot\widetilde{\Pi}_{0}\cdot\sigma_0+\sigma^{-1}\cdot\widetilde{\Pi}\cdot\sigma
\end{equation}
is more naturally expressed in terms of the projectors $\widetilde{\Pi}_{0}\colon \mathbb{R}\to\mathbb{R}$ and $\widetilde{\Pi}\colon \mathbb{C}^{N-1}\to\mathbb{C}^{N-1}$. The projectors $\Pi_B$ and $\widetilde{\Pi}$ act componentwise and the action of all three projectors on components $\rho_i,\tilde{\rho}_0\in\mathbb{R}$ and $\tilde{\rho}_j\in\mathbb{C}$ takes a similar form:
\begin{eqnarray}
\Pi_B(\rho_i) &:=& \left\{\begin{array}{ll}
1/2(\rho_i/|\rho_i|)  & \mbox{if $\rho_i\ne 0$},\cr
1/2  & \mbox{otherwise.}\end{array}
\right.\quad(0\leq i\leq N-1)\\
\widetilde{\Pi}_0(\tilde{\rho}_0) &:=& \left\{\begin{array}{ll}
\sqrt{\sigma_0(\alpha)}(\tilde{\rho}_0/|\tilde{\rho}_0|)  & \mbox{if $\tilde{\rho}_0\ne 0$},\cr
\sqrt{\sigma_0(\alpha)}  & \mbox{otherwise.}\end{array}
\right.\\
\widetilde{\Pi}(\tilde{\rho}_j) &:=& \left\{\begin{array}{ll}
\sqrt{\sigma_j(\alpha)}(\tilde{\rho}_j/|\tilde{\rho}_j|)  & \mbox{if $\tilde{\rho}_j\ne 0$},\cr
\sqrt{\sigma_j(\alpha)}  & \mbox{otherwise.}\end{array}
\right.\quad(1\leq j\leq N-1)
\end{eqnarray}
That all three are distance minimizing is immediately clear given the two ways (\ref{quadForm1}, \ref{quadForm3}) of expressing the Euclidean norm; the definitions for the exceptional cases ($\rho_i=0$, etc.) are arbitrary but apply to sets of measure zero and therefore never arise in actual computations.

\begin{figure}
\begin{center}
\scalebox{1.}{\includegraphics{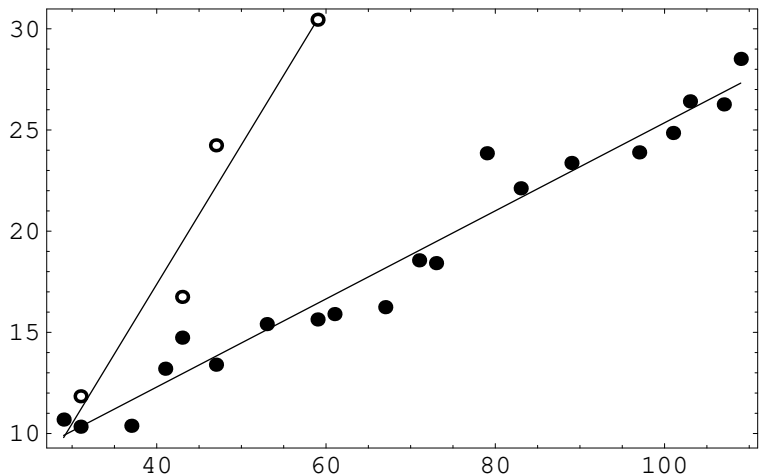}}
\end{center}
\begin{center}
\parbox{5in}{
\caption{Complexity of the difference map algorithm for two sets of bit retrieval instances. Plotted vertically is $\log_2{(I_0)}$, where $I_0$ is the mean number of iterations performed by the algorithm. Instances fall in the range $29\leq N\leq 109$ (horizontal axis) and include $\pi$-sequences (solid circles) and Hadamard sequences (open circles).}}
\end{center}
\end{figure}

The difference map with parameter $\beta=0.7$ (see Appendix) found solutions for bit retrieval instances significantly faster than either of the other algorithms (Table 2).
Figure 1 shows results for the $\pi$-sequence instances in the range $29\leq N\leq 109$ and the significantly more difficult Hadamard sequences for $N=31, 43, 47$ and 59.
Several runs were performed for each instance in order to reliably obtain the mean number of iterations $I_0$ required by the algorithm to find the solution. From the overall linear variation of $\log_2{(I_0)}$ with $N$ of the $\pi$-sequence instances, one obtains the estimate $2^{c N}$ for the average-case complexity, with $c\approx 0.22$. The complexity is dominated by the exponential number of iterations performed, since the time required per iteration grows only as $O(N\log{N})$ (from FFT computations). The Hadamard sequences were selected for study because they saturate the norm bound (theorem \ref{normBoundThm}). For these instances the complexity of the algorithm follows a distinctly steeper exponential growth, with $c\approx 0.69$.

\section{Public key signature}\label{signature}

The economy of hiding binary sequences within their autocorrelation almost rivals that of the RSA scheme of hiding a pair of large primes within their product [RSA]. As for the task of retrieving binary sequences from their autocorrelation, the survey of algorithms in the previous section lends some evidence to the possibility that bit retrieval may be even harder than factoring large integers. These two considerations combined, economy and intractability, provide motivation to design cryptographic systems based on the one-way nature of the autocorrelation operation. Below we propose a digital signature where private and public keys are related by this one-way function. In its broadest description this scheme belongs to the class of cryptographic systems based on lattices (see [MG]), a notable example being the NTRU system [NTRU] whose lattices, as here, are ideals of the ring $Z$. The characteristic of the new scheme that represents a departure from other lattice-based systems, including NTRU, is the simplicity of the relationship between private and public keys. In that the latter can be viewed as the product in an algebraic number field, the RSA relationship between private and public keys provides a natural point of comparison. On the other hand, by using the degree of the number field ($N-1$) as the security parameter, and in particular not having the benefit of a Euclidean division algorithm, the new scheme enters largely unexplored territory. 

A brief description of the scheme developed below begins with Alice, who wishes to apply her signature to a piece of data. We consider two closely related situations: (1) Alice signs a general digital document by \textit{attaching} her signature, and (2) Alice signs data that may even be analog in nature by \textit{modifying} it irreversibly. The term \textit{watermark} will be used when referring to case (2). In both cases the input to the signing operation is an element $\rho\in R$. The watermarking situation is the most straightforward, where $\rho$ is simply a set of $N$ samples of say an audio signal or grayscale image. We assume the individual samples are measured to sufficient resolution such that when rescaled to unit resolution the corresponding elements $\rho\in Z$ have a large range, say $0<[\rho]_i <M$ with $M=2^8$, for example. In the more general situation (1), we assume that the element $\rho\in Z$ is the output of a public message digest (one-way hash function), applied to the digital document.

Alice's private key is a secret binary integer $\beta\in O$ that defines a map $S_\beta\colon R\to Z$ which sends the input $\rho$ to an element $\rho_\beta\in Z$ with the property $\Psi(\rho_\beta)\in\beta\, O$. In essence, the signing operation corresponds to quantization of the ``cyclotomic content" of $\rho$ on a secret principal ideal. A key property of the signing map is the guarantee $\|\rho-\rho_\beta\|<\Delta$, where $\Delta$ is a parameter. In the watermarking scenario this is clearly important if the signed data is to serve as a substitute for the original. More significantly, particularly when signing a message digest for which fidelity is not an issue, the smallness of $\Delta$ provides security against forgeries.

By signing the data Alice hopes to be able to assert her authorship when challenged, for example, by Bob. Moreover, Bob may independently have an interest in establishing the authenticity of data attributed to Alice. Both needs are met if Alice publishes the autocorrelation of her private key, $\alpha=\beta\overline{\beta}$. To verify authorship or authenticity, Bob must check two things. First, he computes the autocorrelation of the data in $O$, $\Psi(\rho_\beta\overline{\rho}_\beta)$, and checks for divisibility by Alice's public key $\alpha$. If $\alpha$ does not divide $\Psi(\rho_\beta\overline{\rho}_\beta)$, then Bob concludes the data is not quantized on Alice's secret ideal $\beta\,O$. Second, in the message digest scenario, Bob applies the public hash function to the document to obtain $\rho$ and checks that $\|\rho-\rho_\beta\|<\Delta$. If the inequality is violated Bob concludes that the signature was forged. In the watermarking scenario, where Bob does not have access to the original $\rho$, the violation of this inequality manifests itself in a signal, image, etc. that is so distorted to be immediately suspect. 

The security of this scheme rests on two assumptions: (1) extracting Alice's private key from her public key, or bit retrieval, is computationally infeasible, and (2) without access to Alice's private key it is infeasible to compute good quantizers for her secret ideal. Attacks which test these assumptions will be refereed to as ``direct" and ``counterfeiting", respectively.

\subsection{Key generation}

From the empirical complexity estimate $2^{c N}$, $c\approx 0.22$, for the fastest known algorithm, it appears that bit retrieval becomes effectively infeasible for relatively modest values of $N$, say $N>250$. Once $N$ is fixed, the success of bit retrieval by the difference map can be further diminished by increasing the norm of the private key $\beta$, as implied by the observed correlation between the latter and the average number of iterations performed by the algorithm (Fig. 1). Since the norm can be calculated efficiently, a practical method for optimizing the key is to simply generate a large number of binary integers using a pseudo-random number generator and select the one with the largest norm.

\subsection{Signing}

The process of signing an element $\rho\in R$ (data, message digest) is accomplished by the map $S_\beta\colon R\to Z$ defined by
\begin{equation}\label{signingMap}
S_\beta(\rho):=\lceil\sigma^{-1}\cdot\sigma\left(Q_\beta(\rho)\right)+\pi_0(\rho)\rfloor\; ,
\end{equation}
where $Q_\beta\colon R\to \beta\, O$ is the quantizing map that requires the private key $\beta$, and $\lceil\;\rfloor$ rounds each component in the standard basis to the nearest integer. Since $Q_\beta(\rho)\in O$, we have $\sigma^{-1}\cdot\sigma\left(Q_\beta(\rho)\right)=\alpha+q\, \Phi_N$ for some $\alpha\in Z$ and $q\in\mathbb{Q}$. Moreover, since $\pi_0(\rho)=r\,\Phi_N$ for some $r\in\mathbb{R}$, all components acted upon by the rounding operation have the same fractional part and we have
\begin{equation}\label{sign1}
S_\beta(\rho)=\sigma^{-1}\cdot\sigma\left(Q_\beta(\rho)\right)+\pi_0(\rho)+\epsilon\,\Phi_N\; ,
\end{equation}
where $|\epsilon|<\frac{1}{2}$. From (\ref{sign1}) we infer that $\Psi(S_\beta(\rho))=Q_\beta(\rho)$ and $\pi_0(S_\beta(\rho)-\rho)=\epsilon\,\Phi_N$, showing that $S_\beta$ preserves the cyclotomic ``codeword" $Q_\beta(\rho)$ and the embedding in $Z$ achieves the minimum distance when projected onto the ideal $\mathbb{R}\,\Phi_N$.

The quantizing map $Q_\beta$ seeks to find the element of the ideal $\beta\, O$ that minimizes the Euclidean distance to $\rho$ in the orthogonal complement of $\mathbb{R}\,\Phi_N$, the space $R_\perp\cong R/\langle\Phi_N\rangle$. Since this \textit{closest vector} problem is hard for the arbitrary ideals (lattices) specified by $\beta$, we use an approximate but computationally efficient form for $Q_\beta$. For arbitrary $\rho\in R$, define
\begin{equation}\label{Qbeta}
Q_\beta(\rho):=\beta\,Q_O\left(\sigma^{-1}\cdot(\sigma(\rho)/\sigma(\beta))\right)\; ,
\end{equation}
where the division sign denotes componentwise division and $Q_O$ is the quantizer $R_\perp\to O$ for the norm (\ref{Onorm}). For $\beta\ne 0$ this map is well defined since the complex numbers $\sigma_j(\beta)$ will all be nonzero.

The problem of computing $Q_O(\gamma)$ for $\gamma\in R_\perp$ is equivalent to vector quantization for the dual of the root lattice $A_{N-1}$ and is treated by Conway and Sloane [CS]. In the following we describe the algorithm given by Scheidler and Williams [SW] in the context of Euclidean division algorithms for cyclotomic fields. We first obtain $\lfloor\gamma\rfloor\in Z$ by taking the floor of each component in the standard basis. The fractional parts of the components are then sorted to obtain a permutation $\{p_0\,\ldots\, p_{N-1}\}$ of $\{0\,\ldots\, N-1\}$ such that if $\gamma-\lfloor\gamma\rfloor=\sum_{i=0}^{N-1}\epsilon_i\, x^{p_i}$, then $\epsilon_0\leq\epsilon_1\cdots\leq\epsilon_{N-1}$. Using this permutation we recursively generate the sequence $\gamma_0\,\ldots\,\gamma_{N-1}$, where $\gamma_0=\lfloor\gamma\rfloor$ and $\gamma_{i+1}=\gamma_i+x^{p_i}$. The quantizer is then given by $Q_O(\gamma)=\Psi(\gamma_i)$, where $i$ identifies the element of the sequence that minimizes $\|\gamma-\gamma_i\|_\perp$. From the geometry of the fundamental domain $D\subset R_\perp$ of $O$ (see [CS], [L]) one obtains the following bound on the quantization error:
\begin{equation}
\|\gamma-Q_O(\gamma)\|_\perp\leq\frac{N^2-1}{12\,N}\; .
\end{equation}
The \textit{mean-squared quantization error} $\Delta_O$ is defined as the expectation value of $\|\gamma-Q_O(\gamma)\|_\perp$ when $\gamma$ is uniformly distributed over a region in $R_\perp$ that is large enough that edge effects can be neglected, or equivalently, where $\gamma$ is uniformly distributed over $D$. A formula for $\Delta_O$, useful for small $N$, is given in [CS].

When $N$ is large a good alternative to the quantizer $Q_O$ is the simpler map $Q_Z\colon\gamma\mapsto\Psi(\lceil\gamma\rfloor)$. For uniformly distributed data one can show [E2] that the improvement in the quantization error, of $Q_O$ over $Q_Z$, is almost always negligible as $N\to \infty$, a fact that also implies the asymptotic limit $\Delta_O\sim \frac{N}{12}$. The approximate quantizer $Q_Z$ can be computed somewhat faster than $Q_O$.

A quantitative measure of the fidelity of the signed data is the evaluation of the mean-squared quantization error $\Delta_\beta$ of the map $S_\beta$:
\begin{thm}
For the signing map specified by (\ref{signingMap}) and uniformly distributed data $\rho\in R$,
\begin{eqnarray}
\Delta_\beta &:=& \mathrm{E}\left(\|S_\beta(\rho)-\rho\|\right)\\
&=& N\left( \frac{\Delta_O}{N-1}\, \|\beta\|_\perp+\frac{1}{12}\right)\label{errorProof0}
\end{eqnarray}
\end{thm}
\begin{proof}
The calculation combines the two forms of quantization error already discussed. The arguments leading to (\ref{sign1}) show that the real number $\epsilon$ defined by
\begin{equation}\label{errorProof1}
S_\beta(\rho)-\rho=\sigma^{-1}\cdot\left(\sigma(Q_\beta(\rho)-\sigma(\rho)\right)+\epsilon\,\Phi_N
\end{equation}
is uniformly distributed in the interval $(-\frac{1}{2},\frac{1}{2})$ when $\rho$ is uniformly distributed in $R$. For quantization in $R_\perp$ by $Q_O$ we have the statement that $\delta\in R_\perp$ defined in
\begin{equation}\label{errorProof2}
Q_\beta(\rho)=\beta\left(\sigma^{-1}\cdot(\sigma(\rho)/\sigma(\beta))+\delta\right)
\end{equation}
is uniformly distributed in the fundamental region $D$ of $O$. Moreover, the distributions of $\epsilon$ and $\delta$ are clearly independent. From (\ref{errorProof2}) we have
\begin{equation}
\sigma(Q_\beta)=\sigma(\rho)+\sigma(\beta)\sigma(\delta)\; ,
\end{equation}
with the result that (\ref{errorProof1}) may be rewritten as
\begin{equation}\label{errorProof3}
S_\beta(\rho)-\rho=\sigma^{-1}\cdot(\sigma(\beta)\sigma(\delta))+\epsilon\,\Phi_N\; .
\end{equation}
Taking the norm of (\ref{errorProof3}) we have
\begin{eqnarray}
\|S_\beta(\rho)-\rho\| &=& \|\sigma^{-1}\cdot(\sigma(\beta)\sigma(\delta))\|_\perp+\epsilon^2\, N\\
&=& \frac{1}{N}\left(\sigma(\beta)\sigma(\delta)\right)\cdot \left(\overline{\sigma}(\beta)\overline{\sigma}(\delta)\right)+\epsilon^2\, N\; .\label{errorProof4}
\end{eqnarray}
What remains is taking the expectation values $\mathrm{E}(\epsilon^2)=\frac{1}{12}$ and for $1\leq j\leq N-1$,
\begin{eqnarray}
\mathrm{E}\left(\sigma_j(\delta)\overline{\sigma}_j(\delta)\right) &=& \frac{1}{N-1}\,\mathrm{E}\left(\sigma(\delta)\cdot\overline{\sigma}(\delta)\right)\label{errorProof5}\\
&=& \frac{N}{N-1}\,\mathrm{E}(\|\delta\|_\perp)\\
&=& \frac{N}{N-1}\,\Delta_O\; ,
\end{eqnarray}
since the left side of (\ref{errorProof5}) is clearly independent of $j$. After applying these averages to (\ref{errorProof4}) we obtain the result (\ref{errorProof0}) for the mean-squared quantization error:
\begin{eqnarray}
\Delta_\beta &=& \frac{1}{N}\sum_{j=1}^{N-1}\sigma_j(\beta)\overline{\sigma}_j(\beta)\,\mathrm{E}\left(\sigma_j(\delta)\overline{\sigma}_j(\delta)\right)+\mathrm{E}(\epsilon^2) \,N\\
&=& \frac{\Delta_O}{N-1}\sum_{j=1}^{N-1}\sigma_j(\beta)\overline{\sigma}_j(\beta)+\frac{N}{12}\; .
\end{eqnarray}
\end{proof}

The most direct way of assessing the fidelity of a watermark is by comparing the \textit{root-mean-squared quantization error per component} for uniformly distributed data,
\begin{equation}\label{rmsError1}
\delta_{\mathrm{rms}}:=\sqrt{\frac{\Delta_\beta}{N}}\; ,
\end{equation}
with the range of values in the data. We are primarily interested in $\delta_{\mathrm{rms}}$ when $\beta$ is a random binary key and $N$ is large. If $\beta_R$ is the corresponding binary key in $R$, then $\|\beta\|_\perp=\|\beta_R\|-\sigma_0(\beta_R)^2/N\sim N/4$, since $\sigma_0(\beta_R)=O(\sqrt{N})$. Combining this with (\ref{errorProof0}) and $\Delta_O\sim \frac{N}{12}$, we obtain
\begin{equation}\label{rmsError2}
\delta_{\mathrm{rms}}\sim\sqrt{\frac{N}{48}}\quad (N\to\infty)\; .
\end{equation}
In the image watermarking application of section \ref{watermarking}, for example, the elements of data are blocks of $N=379$ pixels, and the range of each component (pixel) is an 8-bit integer. Signing an image with the map $S_\beta$ thus modifies each pixel ($\pm$) by $\delta_{\mathrm{rms}}\approx 2.8$, or about 1\% of its range.

Associated with the application of a watermark is a loss of information that can be used as a means of normalization when comparing with other schemes. The map $S_\beta$ is an example of a lattice quantizer, for which the lost information content corresponds to the volume $V$ of the region in $R$ that maps to any particular ``codeword" in $Z$. Since this region comprises the product of a unit interval in $\mathbb{R}\Phi_N$ with a fundamental region of $\beta\,O$ in $R_\perp$, we have $V=\mathcal{N}(\beta)$. The standard normalization applied to the root-mean-square quantization error per component is the following [CS]:
\begin{eqnarray}
G &:=& \frac{\delta_{\mathrm{rms}}^2}{V^{2/N}}\\
&\sim& \frac{e^\gamma}{12}\approx 0.148423\ldots\quad (N\to\infty)\; ,\label{G}
\end{eqnarray}
where (\ref{G}) was obtained using (\ref{rmsError2}) and (\ref{E(S)}) for random binary keys $\beta$. 
When the input to the signing operation is already a digital document, this value can be compared with Wong's watermarking scheme [W]. In Wong's scheme the least significant bit of each element of a block of data is replaced by the output of a one-way hash function applied to the block. The parameters for Wong's watermark are thus $\delta_{\mathrm{rms}}=\frac{1}{2}$, $V=2^N$, giving the slightly better value $G=\frac{1}{8}$. On the other hand, for analog data Wong's watermark can only be applied after a digital encoding step has made its own contribution to the net quantization error. For large $N$, Zador's analysis of random quantizers [Za] gives the bound $G>1/(2\pi e)$. 

A noteworthy property of the signing operation, as well as the verification step (below), is that it can be efficiently implemented without the need for arbitrary precision arithmetic: a finite precision general purpose FFT can perform all the necessary ring multiplications and divisions in a time that grows as $N\log{N}$. Assuming that the Fourier transform coefficients of the key, $\sigma(\beta)$, are computed only once during the signing of many data items, a total of four FFTs are performed in the computation of $S_\beta(\rho)$ for each $\rho$. Since all the other parts of the computation (quantizing with $Q_Z$, etc.) only involve $O(N)$ arithmetic operations, the overall complexity of signing is nearly linear in the size of the data, $O(N\log{N})$. Verification is the stronger test of the finite precision arithmetic in that autocorrelations are involved. Tests with 12-bit data showed that standard double precision arithmetic was adequate for $N<1000$.

\subsection{Verification}

To verify that a digital document has been signed by Alice, Bob makes use of four things: the message digest $\rho\in R$ resulting from the application of a public one-way hash function to the document, Alice's signed modification $\rho_\beta=S_\beta(\rho)\in Z$, Alice's public key $\alpha\in \hat{O}$, and the fidelity parameter $\Delta$. He first applies the verification map $V_\alpha\colon Z\to R_\perp$
\begin{equation}
V_\alpha(\rho_\beta):=\sigma^{-1}\cdot\left(\sigma(\rho_\beta\overline{\rho}_\beta)/\sigma(\alpha)\right)\; ,
\end{equation}
and checks whether $V_\alpha(\rho_\beta)\in O$. Recall that if $\rho_\beta$ is quantized with Alice's private key $\beta$, then $\sigma(\rho_\beta)=\sigma(\beta\gamma)$ for some $\gamma\in O$. Since $\alpha=\beta\overline{\beta}$, Bob computes
\begin{equation}
V_\alpha(\rho_\beta)=\sigma^{-1}\cdot\left(\sigma(\beta\gamma\overline{\beta\gamma})/\sigma(\beta\overline{\beta})\right)=\sigma^{-1}\cdot\sigma(\gamma\overline{\gamma})
\end{equation}
and concludes that $V_\alpha(\rho_\beta)\in O$. When unsuccessful, $V_\alpha(\rho_\beta)$ is a non-integer in the cyclotomic field $\mathbb{Q}[\zeta]$, that is, not all components in the standard basis will be integers.

A fast, finite precision arithmetic implementation of this first part of the verification requires two FFTs, not counting $\sigma(\alpha)$, which is computed once in the course of verifying a large stream of data. With the first FFT Bob computes $\sigma(\rho_\beta)$; he then squares the modulus, divides by $\sigma(\alpha)$, and applies the inverse FFT to the result. To check for membership in $O$, he obtains the fractional parts of the components in the standard basis and compares these with zero, making allowance for the finite precision in the calculation.

To complete the verification Bob checks that $\|\rho-\rho_\beta\|<\Delta$. The parameter $\Delta$ is chosen to guard against forgeries. As discussed below, there is a significant gap between the range of distances $\|\rho-\rho_\beta\|$ realized by Alice's quantizers $\rho_\beta$ and quantizers that can be computed by a forger. This gap grows with $N$ so that $\Delta$ need not be specified precisely when $N$ is large. In watermarking applications the last step of the verification cannot be performed because the original $\rho$ is not available. Instead, the poorness of the forger's quantizers have the effect of introducing so much noise to the signal or image that the authenticity of the signature is immediately called into question (see section \ref{watermarking}).

\section{Security}\label{security}

Eve has at least two ways of undermining this signature scheme: she can attempt to determine Alice's private key $\beta$ from the publicly available data, or she can sign data with a substitute for Alice's key and hope that nobody notices. It appears that both forms of attack, respectively direct and counterfeiting, become prohibitively difficult for reasonable values of $N$. 

\subsection{Direct attack}

Since Eve has access to Alice's public key $\alpha=\beta\overline{\beta}$, as well as multiple signed data elements, $\rho_1=\beta\gamma_1$, \mbox{$\rho_2=\beta\gamma_2,\ldots$}, it is fortunate (for Alice)
that Euclidean algorithms cease to exist beyond $N=19$ [MM] that Eve might use to extract the common divisor $\beta$. An alternative approach for solving these instances of problems \textbf{B$_2$} and \textbf{B$_3$} is to use the algorithms of algebraic number theory, as illustrated in section \ref{algNumAlgorithm}. However, neither this approach nor the integer programming method for solving \textbf{B$_2^\prime$} was found to be competitive with the phase retrieval algorithm. The time complexity of the latter was investigated in section \ref{diffMapAlgorithm} and appears to be exponential in $N$. A direct attack is thus infeasible with the currently known algorithms.

\subsection{Counterfeiting}

Since the verification challenge for this signature scheme tests for membership in the ideal $\beta\,O$, and the inclusion $\beta\gamma\,O\subset \beta\,O$ holds for arbitrary $\gamma\in O$, data signed with any nonzero multiple of Alice's private key, say $\beta^\prime=\beta\gamma$, will also satisfy the challenge. Such counterfeit keys are publicly available, from Alice's public key $\alpha=\beta\overline{\beta}$, to the numerous elements of data Alice herself has signed: $\rho_1=\beta\gamma_1$, \mbox{$\rho_2=\beta\gamma_2,\ldots$}\,. What makes these options for counterfeiting Alice's signature generally unacceptable to Eve is that the corresponding quantization errors will be large. The derivation of the root-mean-squared quantization error per component (\ref{rmsError1}) is valid for arbitrary keys $\beta^\prime$ (not necessarily binary) and can be approximated for large $N$ by
\begin{equation}
\delta_{\mathrm{rms}}\sim\sqrt{\frac{\|\beta^\prime\|_\perp}{12}}\; .
\end{equation}
Now if $\beta^\prime=\beta$ is a genuine (binary) private key, then the expectation value
\begin{equation}
\mathrm{E}(\|\beta\|_\perp)\sim \frac{N}{4}\; ,
\end{equation}
assuming a uniform distribution on the binary keys, gives the estimate $\delta_{\mathrm{rms}}\sim\sqrt{N/48}$ obtained previously in (\ref{rmsError2}). If instead $\beta^\prime=\beta\gamma$, then
\begin{equation}
\mathrm{E}(\|\beta\gamma\|_\perp)\sim \frac{N}{4}\|\gamma\|_\perp\; ,
\end{equation}
where the expectation value is again computed (details omitted) with respect to the uniform distribution on binary $\beta$. The counterfeit key thus increases $\delta_{\mathrm{rms}}$ by a factor of order $\sqrt{\|\gamma\|_\perp}$.
If instead Eve chooses to sign with Alice's public key, $\beta^\prime=\beta\overline{\beta}$, then the expectation value (over uniformly distributed binary $\beta$)
\begin{equation}
\mathrm{E}(\|\beta\overline{\beta}\|_\perp)\sim \frac{N^2}{8}
\end{equation}
shows that $\delta_{\mathrm{rms}}$ would increase by $\sqrt{N/2}$ over its value when signing with the private key.

The discussion above suggests making a minor modification to the quantization map (\ref{Qbeta}) that ensures the outcomes $Q_\beta(\rho_\perp)=\beta\gamma$ have factors $\gamma$ with some minimum Euclidean norm that grows with $N$. It was already argued that fidelity is not significantly sacrificed when the quantizer $Q_O$ is replaced by $Q_Z$, and in fact this can be generalized to include the quantizer [E2]
\begin{equation}
Q_{\widetilde{Z}}\colon \rho_\perp\mapsto \Psi(\lceil\rho_\perp+r\,\Phi_N\rfloor)
\end{equation}
for arbitrary $r\in \mathbb{R}$. The choice $r=\frac{1}{2}$ has a clear advantage when signing a block of data $\rho$ where the components are nearly equal, as in watermarking parts of an image with small contrast. The input $\rho_\perp=\sigma^{-1}\cdot\sigma\left(\sigma(\rho)/\sigma(\beta)\right)$ to $Q_{\widetilde{Z}}$ is then a random vector with small components distributed around zero, for which $Q_{\widetilde{Z}}$ with $r=\frac{1}{2}$ produces a random binary integer as output. Quantizing with $Q_{\widetilde{Z}}$ will thus almost always produce factors $\gamma$ with $\|\gamma\|_\perp>N/4$. In the rare event that this is not true, the signer (Alice) can artificially amplify the contrast (by rescaling $\rho_\perp$) until this condition is met.

Eve is also severely limited in how much she can reduce her quantization error through the use of a better quantization algorithm. Since the dimensionless mean-squared quantization error is always greater than Zador's bound $G>1/(2\pi e)$ [Za], and Alice's quantizer $S_\beta$ has $G=e^\gamma/12$, Eve can at most hope to reduce $\delta_{\mathrm{rms}}$ by the constant factor $\sqrt{6/(\pi e^{1+\gamma})}\approx 0.628$.

Eve can mount a different counterfeiting attack by attempting to solve problem \textbf{B$_3$}. Suppose $\gamma_1$ and $\gamma_2$ are two random elements of $O$, say with bounded components. For large $N$ it will almost always be true that $\gamma_1\,O+\gamma_2\,O=O$. Since Eve has access to several products $\beta\gamma_1$, $\beta\gamma_2,\ldots$ (signed data and public key), she can in principle construct the ideal generated by Alice's private key from the fact $\beta\gamma_1\,O+\beta\gamma_2\,O=\beta\,O$. In computational terms this corresponds to taking the union of the lattice generators of the two ideals and applying some form of lattice basis reduction in order to be able to recognize $\beta$. For counterfeiting purposes, however, Eve does not have to succeed in finding $\beta$: rather, she will be satisfied with any element of $\beta\,O$ having a small Euclidean norm. Since this is exactly the kind of problem for which LLL basis reduction has proven to be effective, the following experiment was performed.

For each $N$ in the experiment, twenty ``LLL attacks" were performed. The data for each attack was generated from three random binary integers: $\beta$ (the private key), $\beta_1$ and $\beta_2$. Available to Eve are the pair, $\rho_1=\beta\beta_1$ and $\rho_2=\beta\beta_2$, representing two signed elements of data with small Euclidean norm, say, or one data item and the public key. The lattice basis $\Gamma$ for $\rho_1\,O+\rho_2\,O$ was constructed from $\Gamma_1$ and $\Gamma_2$, where
\begin{equation}
\Gamma_k=\{N\sigma^{-1}\cdot\sigma(\rho_k \zeta^i)\colon 1\leq i\leq N-1\}\quad (k=1, 2)\; .
\end{equation}
The scaling factor $N$ produces an integral basis $\Gamma$ for a lattice in $R_\perp$ to which basis reduction can be applied. From the construction of $\Gamma$ it will almost always be true that there exists a reduced basis $\Gamma^\prime$ where all the generators are binary vectors ($\beta\zeta^i$, $1\leq i\leq N-1$) multiplied by $N$. The minimum Euclidean norm achieved for this reduction (after division by the scaling factor $N$) is therefore $\|\beta\|_\perp\sim N/4$. Computing $\Gamma^\prime$ from $\Gamma$ is difficult, and we limit ourself to the reduced basis $\Gamma_\mathrm{LLL}$ obtained by the LLL algorithm [LLL]. If the basis element of minimal norm, $\gamma_\mathrm{min}\in \Gamma_\mathrm{LLL}$, has an acceptably small norm, Eve can use it as a counterfeit key. As a figure of merit, the output of the experiment was the smallest value of the ratio $r=\|\gamma_\mathrm{min}\|_\perp/(N/4)$ achieved for all twenty attacks. Values $r\approx 1$ indicate a successful attack, that is, where data signed with $\gamma_\mathrm{min}$ would not be noticeably more distorted than data signed with Alice's private key. Unsuccessful attacks have $r>1$ , where larger values result in signed data that is more easily recognized as bearing a counterfeit signature.

Figure 2 shows a plot of $r$ for the range $23\leq N\leq 97$. For $N<50$ the LLL attack is successful, providing Eve with a key in a reasonable time with which she can sign data that would be verified as Alice's. Beyond $N\approx 50$ the ratio $r$ achieved by the LLL attack increases sharply to values where the computed key is not useable. Interestingly, for $N\ge 89$ it appears that LLL basis reduction is even counterproductive, the resulting $\gamma_\mathrm{min}$ having a norm that exceeds the norms $\|\rho_1\|_\perp\approx\|\rho_2\|_\perp\sim (N/4)^2$ of the starting basis elements (shown as the line with slope $1/4$ in Fig. 2).

\begin{figure}
\begin{center}
\scalebox{1.}{\includegraphics{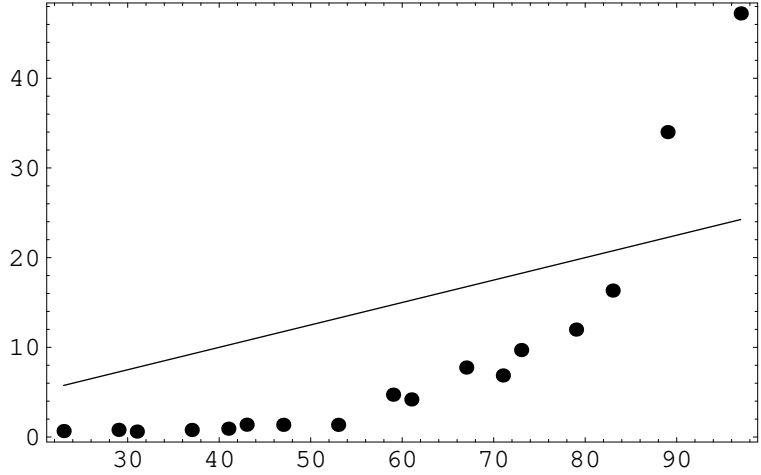}}
\end{center}
\begin{center}
\parbox{5in}{
\caption{Failure of LLL basis reduction to find a suitable counterfeit key when $N$ (horizontal axis) is large. Each data point represents the smallest norm basis element $\gamma_\mathrm{min}$ found by the LLL algorithm out of twenty trials. The vertical axis is the ratio $r=\|\gamma_\mathrm{min}\|_\perp/(N/4)$, or the excess norm over the reduction corresponding to the discovery of the private key.}}
\end{center}
\end{figure}

\section{Image watermarking}\label{watermarking}

The signature scheme proposed in the previous section, when applied to image watermarking, illustrates the role of noise in the detection of forgeries. We recall that increasing the value of the security parameter $N$ serves two purposes: (1) the corresponding bit retrieval problem, of extracting the private key from the public key or signed data elements, becomes harder, and (2) the quality of quantization with counterfeit keys becomes increasingly poor. Here we focus entirely on the second point. 

The creation of forgeries in the present context is known in the watermarking literature as a \textit{vector quantization attack} [HM]. Wong [W] introduced the watermarking scheme where Alice modifies each block of pixels in their least significant bits by the output of a message digest applied to the block. The forger, Eve, is then limited to building her images out of exact copies of blocks that have already appeared in images signed by Alice. 
The set of available image quantizers --- blocks bearing a valid signature --- in the present scheme is considerably larger, being any elements of the lattice specified by Alice's private key. 

There are numerous practical issues that our discussion omits, such as the method of partitioning the image into data blocks [Ce]. We are only interested in watermarks that are both \textit{invisible} and \textit{fragile}. The latter term refers to the property that changes in the value of even one pixel will cause a failure in the verification and facilitate the localization of tampering.

\begin{figure}
\begin{center}
\scalebox{1.6}{\includegraphics{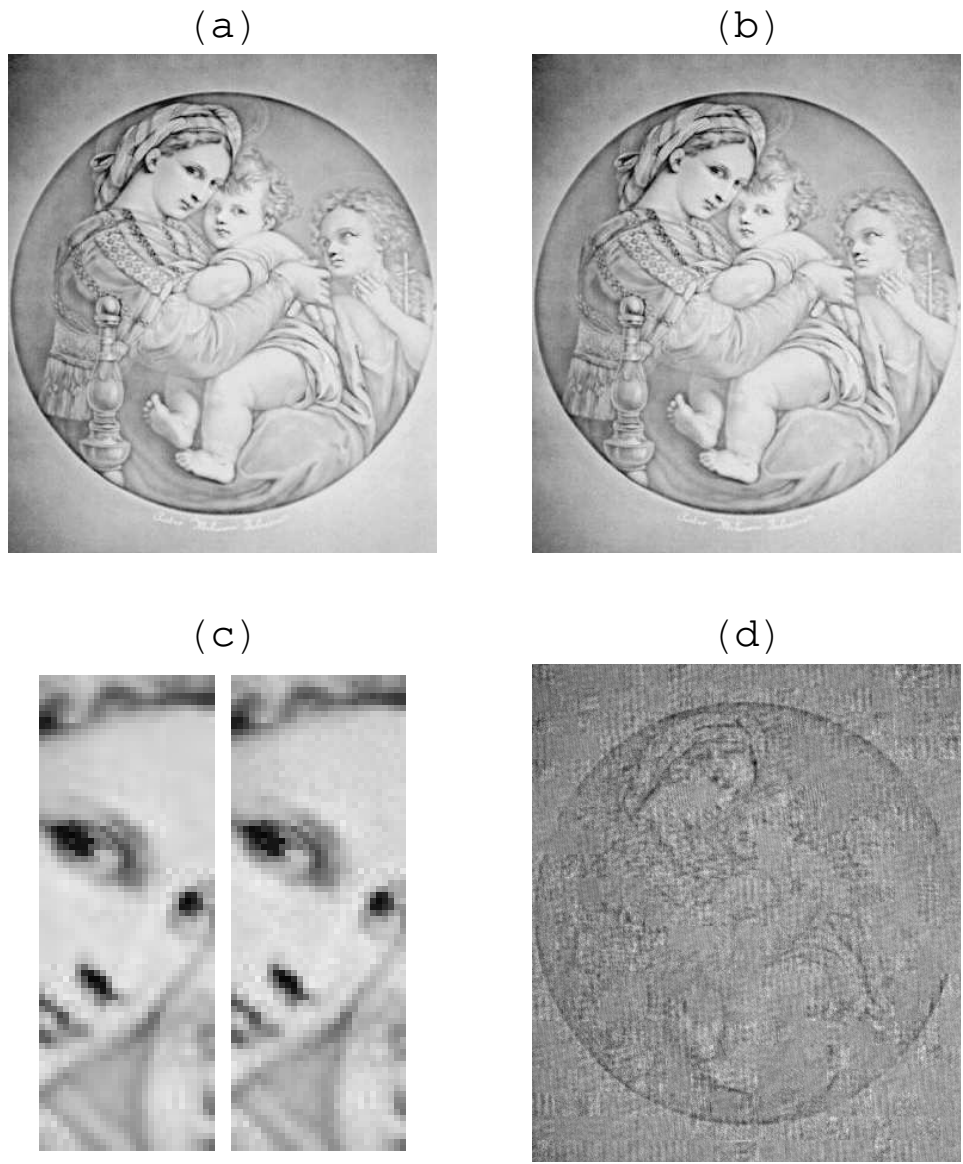}}
\end{center}
\begin{center}
\parbox{5in}{
\caption{Image watermarking application of the digital signature. (a) TIFF image of a paper watermark by Pietro Miliani Fabriano. (b) Modification of (a), signed with a binary key. (c) Details of original (left) and signed (right) images. (d) Noisy image produced by signing (a) with a counterfeit key.}}
\end{center}
\end{figure}

Figure 3 shows the result of applying a digital signature, of the type described in section \ref{signature}, to a $361\times 420$ pixel grayscale image. The pixels of the image were first partitioned into $19\times 20$ rectangular blocks, where the dimensions were chosen so that the total number of pixels per block is one greater than a prime, in this case $N=379$. The extra pixel was left unchanged by the signing operation. To ensure that the final signed image has 8-bit integer pixels, a global scaling and shift was applied to all the pixel values of the original. Since signing typically modifies a component by $\pm\delta_\textrm{rms}\approx 2.8$ (for $N=379$), the parameters of the scaling and shift were adjusted to bring the pixel values of the original into the range $5-250$. Frames (a) and (b) of Figure 3 are TIFF images using, respectively, the original and signed pixels as raster data. The two images are practically indistinguishable, with differences (c) discernible only at artificially high magnification. 

An attempted forgery is shown in (d), where quantization was not performed with Alice's short binary key $\beta$, but a much longer counterfeit key $\beta\gamma$. The latter key was taken from Alice's signed image (b), specifically from the pixel block with smallest Euclidean norm. Blocks with small Euclidean norm arise in those parts of an image where the contrast is small. Recognizing this, image (b) was signed using the quantizer $Q_{\widetilde{Z}}$ (and $r=\frac{1}{2}$) which avoids multipliers $\gamma$ with small norms. In (b) the smallest norm among the 399 blocks, $\|\beta\gamma\|_\perp\approx 4697$, was considerably larger than that of the private key, $\|\beta\|_\perp\approx 95$. The poor quality of the resulting forgery is the result of two mechanisms. First, the amplitude of the noise introduced by signing, or $\delta_\mathrm{rms}$, is increased by the factor $\sqrt{\|\beta\gamma\|_\perp/\|\beta\|_\perp}\approx 7$. Second, the increased value of $\delta_\mathrm{rms}$ requires that the range in the pixel values of the original must first be compressed (by rescaling) in order that the signed values fall in the range 0 -- 255. The second mechanism has the effect of reducing the signal to noise ratio of the signed image to practically zero when the counterfeit key $\beta\gamma$ has a sufficiently large norm.

\section{Acknowledgments}
The author thanks J. Buhler, G. Casella, and S. Chase for helpful discussions. This work was supported by the National Science Foundation under grant ITR-0081775.

\section{Appendix: the difference map}

Let $A$ and $B$ be subsets of an $N$-dimensional Euclidean space $E$. For the application discussed in \mbox{section \ref{diffMapAlgorithm}}, $E$ is the ring $R$. The specification of the sets $A$ and $B$ is computationally easy, while the task of computing the intersection $A\cap B$ is assumed to be difficult. The difference map is defined in terms of projectors $\Pi_A$ and $\Pi_B$, which map an arbitrary $x\in E$ to points in $A$ and $B$ that minimize the Euclidean distances, $\|\Pi_A(x)-x\|$ and $\|\Pi_B(x)-x\|$. Practical algorithms require that both projectors can be computed efficiently for any $x\in E$.

We are interested in solving
\begin{equation}\textbf{find:}
\quad x\in A\cap B\; ,
\end{equation}
or equivalently,
\begin{equation}\label{diffMapProblem}\textbf{find:}
\quad x\in E\quad \textbf{such that}\quad x=\Pi_A(x)=\Pi_B(x)\; .
\end{equation}

The difference map $D\colon E\to E$, defined by [E1]
\begin{equation}\label{diffMap}
D(x):=x+\beta (\Pi_B f_A - \Pi_A f_B)(x)\; ,
\end{equation}
is constructed such that its fixed points are simply related to the solutions of (\ref{diffMapProblem}). Here $\beta\neq 0$ is a real parameter and the maps $f_A, f_B\colon E\to E$ are defined in terms of the basic projectors by
\begin{eqnarray}
f_A &:=& (1+\gamma_A)\Pi_A-\gamma_A\\
f_B &:=& (1+\gamma_B)\Pi_B-\gamma_B\; ,
\end{eqnarray}
where $\gamma_A$ and $\gamma_B$ are two additional real parameters. At a fixed point of $D$, $x^\ast=D(x^\ast)$, we have
\begin{equation}
\Pi_B f_A(x^\ast)=\Pi_A f_B(x^\ast):=x_{\mathrm{sol}}\; ,
\end{equation}
and $x_{\mathrm{sol}}$ evidently solves (\ref{diffMapProblem}) since
\begin{equation}
\Pi_A(x_{\mathrm{sol}})=\Pi_A \Pi_A f_B(x^\ast)=\Pi_A f_B(x^\ast)=x_{\mathrm{sol}}\; ,
\end{equation}
and similarly when acted upon by $\Pi_B$. In general $x_{\mathrm{sol}}\neq x^\ast$, 
and the set of fixed points associated with $x_{\mathrm{sol}}$,
\begin{equation}
(\Pi_A f_B)^{-1}(x_{\mathrm{sol}})\cap (\Pi_B f_A)^{-1}(x_{\mathrm{sol}})\; ,
\end{equation}
is normally a continuum. The set of fixed points is not empty if a solution $x_{\mathrm{sol}}$ exists, since $x_{\mathrm{sol}}$ is itself a fixed point.

The parameters $\gamma_A$ and $\gamma_B$ are chosen to make the fixed points of $D$ attractive. Satisfying this criterion for arbitrary sets $A$ and $B$ and optimizing convergence is in general difficult [E2]. Here we consider two particularly simple examples of the local behavior. First, if the sets $A$ and $B$ are manifolds we can approximate them by affine spaces in the neighborhood of a solution. After translating this solution to the origin, we make the further assumption that the corresponding linear spaces are orthogonal so that the projectors satisfy $\Pi_A\Pi_B = 0$. The difference map then simplifies to
\begin{equation}
D(x)=x-\beta\gamma_A\,\Pi_B(x)+\beta\gamma_B\,\Pi_A(x)\;.
\end{equation}
Optimal convergence to the fixed points of $D$ (the linear space $\ker{\Pi_A}\cap\ker{\Pi_B}$) is obtained when
\begin{equation}\label{gammaOpt}
\gamma_A=-\gamma_B=1/\beta\;,
\end{equation}
although this assumes both $A$ and $B$ have positive dimension. If either space is a point, then $\Pi_A=0$ or $\Pi_B=0$ and, respectively, the optimal $\gamma_B$ or $\gamma_A$ remains undetermined. Since this is the case for the set $B$ in bit retrieval (hypercube), our second example examines this situation. For simplicity we take $N=1$ and sets $A=\mathbb{Z}$ and $B=\{0\}$. The corresponding difference map is given by
\begin{equation}
D(x)=x+\beta\lceil\gamma_B x\rfloor\; ,
\end{equation}
where $\lceil\;\rfloor$ denotes rounding to the nearest integer.
The set of fixed points is the interval $\left(-(2\gamma_B)^{-1}, (2\gamma_B)^{-1}\right)$, where the (trivial) local behavior of $D$ is independent of $\gamma_B$ as already mentioned. However, on a global scale we see that convergence requires that $\beta$ and $\gamma_B$ have opposite signs. In fact, optimal convergence is obtained precisely when $\gamma_B=-1/\beta$, in agreement with (\ref{gammaOpt}). In the absence of a more comprehensive analysis we will adopt the parameter values (\ref{gammaOpt}) suggested by these two examples.

A special case of the difference map first appeared in the context of image reconstruction from Fourier modulus data and an object support constraint. Motivated by ideas from linear control theory, Fienup considered three feedback variants in an iterative scheme, the most successful of which became known as the \textit{hybrid input-output} algorithm [F]. In image reconstruction applications of the difference map, $A$ corresponds to the torus of Fourier modulus constraints, as in bit retrieval, while $B$ is a linear space representing the support of the object in the image. Fienup's formulation made no reference to projectors but coincides exactly with the difference map for the parameter values $\gamma_A=1/\beta$, $\gamma_B=-1$, and $\beta>0$ [E1]. The geometrical representation and generalization of the hybrid input-output iteration, made possible by projectors, was recognized only recently [BCL, E1]. 

When applied to bit retrieval and phase retrieval with atomicity constraints, it is believed [E1] that the dynamics of the difference map is chaotic and strongly mixing. If true, this implies that the starting point of the iterations is largely irrelevant: an initial distribution of starting points very quickly approaches an invariant distribution. This property can be strictly true only in the case of ill-posed instances, when there is no solution. Solutions represent an exceptional situation, a constellation of fixed points ``hidden" within the invariant distribution that the chaotic dynamics is attempting to discover. The strongly mixing hypothesis implies that every iteration is effectively subject to a fixed probability of being within the basin of attraction of a fixed point, after which it quickly converges to an entirely different invariant distribution: the fixed point. Thus the number of iterations $I$ of the method is expected to have the probability distribution
\begin{equation}\label{iterDistribution}
dP(I) = \exp{(-I/I_0)}\,(dI/I_0)\; ,
\end{equation}
where $I_0$ is the mean number.
The method is optimized by minimizing $I_0$ with respect to the parameter $\beta$ for appropriate test problems. Figure 4 compares the histogram of the number of iterations required to solve the bit retrieval instance for the sequence $\pi_{41}$ with the distribution (\ref{iterDistribution}). The data shown represent $10^4$ solution attempts, all successful and differing only in the choice of initial (random) iterate.

\begin{figure}
\begin{center}
\scalebox{1.}{\includegraphics{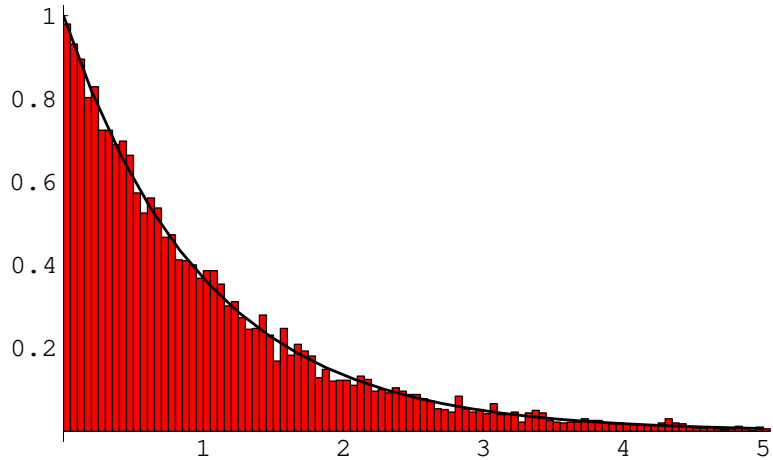}}
\end{center}
\begin{center}
\parbox{5in}{
\caption{Comparison of the distribution of difference map iterations $I$, required to solve the bit retrieval instance $\pi_{41}$, with the exponential distribution predicted by the strongly mixing hypothesis. The units on the abscissa give the ratio $I/I_0$, where $I_0\approx 9623$ is the mean number of iterations.}}
\end{center}
\end{figure}


\begin{thebibliography}{NTRU}
\bibitem[Ba]{Ba} L. D. Baumert, \textit{Cyclic Difference Sets} (Springer-Verlag, Berlin, 1971).
\bibitem[BCL]{BCL} H. H. Bauschke, P. L. Combettes and D. R. Luke, ``Phase retrieval, Gerchberg-Saxton algorithm, and Fienup variants: A view from convex optimization," J. Opt. Soc. Am. A \textbf{19}, 1334-1345 (2002).
\bibitem[Bi]{Bi} P. Billingsley, \textit{Probability and Measure} (John Wiley \& Sons, New York, 1979), p. 310.
\bibitem[Bu]{Bu} J. P. Buhler, private communication.
\bibitem[Ce]{Ce} M. U. Celik, G. Sharma, E. Saber and A. M. Tekalp, ``Hierarchical watermarking
    for secure image authentication with localization,"  IEEE Trans. on Image Proc. \textbf{11} (2002).
\bibitem[Co]{Co} H. Cohen, \textit{A Course in Computational Algebraic Number Theory} (Springer-Verlag, Berlin, 1993).
\bibitem[CS]{CS} J. H. Conway and N. J. A. Sloane, ``Voronoi Regions of Lattices, Second Moments of Polytopes, and Quantization," IEEE Trans. Information Theory, IT-28, 211-226 (1982).
\bibitem[E1]{E1} V. Elser, ``Phase retrieval by iterated projections," J. Opt. Soc. Am. A \textbf{20}, 40-55 (2003).
\bibitem[E2]{E2} V. Elser, unpublished.
\bibitem[E3]{E3} V. Elser, ``Random projections and the optimization of an algorithm for phase retrieval," J. Phys. A: Math. Gen. \textbf{36}, 2995-3007 (2003).
\bibitem[F]{F} J. R. Fienup, ``Phase retrieval algorithms: a comparison,''  Appl. Opt. \textbf{21}, 2758-2769 (1982).
\bibitem[FP]{FP} U. Fincke and M. Pohst, ``Improved methods for calculating vectors of short length in a lattice, including a complexity analysis," Math. Comp. \textbf{44}, 463-471 (1985).
\bibitem[GJ]{GJ} M. R. Garey and D. S. Johnson, \textit{Computers and Intractability, a guide to the theory of NP-completeness} (W. H. Freeman, San Francisco, 1979).
\bibitem[Ha]{Ha} L. Hafer, ``bonsaiG User's Manual," Technical Report SFU-CMPT TR 1999-07, School of Computing Science, Simon Fraser University, Burnaby, B. C., V5A 1S6 (1999).
\bibitem[HM]{HM} M. Holliman and N. Memon, ``Counterfeiting attacks on oblivious block-wise independent invisible watermarking schemes," IEEE Trans. on Image Proc. \textbf{9}, 432-441 (2000).
\bibitem[K]{K} M. Daberkow, C. Fieker, J. Kl\"uners, M. Pohst, K. Roegner and K. Wildanger, ``KANT V4," J. Symbolic Comp. \textbf{24}, 267-283 (1997).
\bibitem[L]{L} H. W. Lenstra, Jr., ``Euclid's algorithm in cyclotomic fields," J. London Math. Soc. \textbf{2}, 457-465 (1975).
\bibitem[LL]{LL} A. K. Lenstra and H.ÊW. Lenstra Jr.,  \textit{The Development of the Number Field Sieve} (Springer-Verlag, Berlin, 1993). 
\bibitem[LLL]{LLL} A. K. Lenstra, H. W. Lenstra Jr. and L. Lov\'asz, ``Factoring polynomials with rational coefficients,"  Math. Ann. \textbf{261}, 515-534 (1982).
\bibitem[LO]{LO} J. C. Lagarias and A. M. Odlyzko, ``Solving low-density subset sum problems," J. ACM \textbf{32}, 229-246 (1985).
\bibitem[MG]{MG} D. Micciancio and S. Goldwasser, \textit{Complexity of Lattice Problems: A Cryptographic Perspective} (Kluwer Academic 2002).
\bibitem[MM]{MM} J. M. Masley and H. L. Montgomery, ``Cyclotomic fields with unique factorization," J. Reine Angew. Math. \textbf{286/287}, 248-256 (1976).
\bibitem[NTRU]{NTRU} J. Hoffstein, J. Pipher and Joseph H. Silverman, ``NTRU: A Ring-Based Public Key Cryptosystem," in Algorithmic Number Theory (ANTS III), Portland, OR, June 1998, J. P. Buhler (ed.), \textit{Lecture Notes in Computer Science} \textbf{1423} (Springer-Verlag, Berlin, 1998) 267-288.
\bibitem[PS]{PS} L. Pauling and M. D. Shappell, Zeits. f. Krist. \textbf{75}, 128 (1930).
\bibitem[R]{R} I. Rankenburg, ``A polynomial-time algorithm for symmetric bit retrieval," unpublished.
\bibitem[RSA]{RSA} R. Rivest, A. Shamir and L. Adleman, ``A method for obtaining digital signatures and public key cryptosystems,"  Comm. ACM \textbf{21}, 120-126 (1978).
\bibitem[SW]{SW} R. Scheidler and H. C. Williams, ``A public-key cryptosystem utilizing cyclotomic fields," Designs, Codes and Cryptography \textbf{6}, 117-131 (1995).
\bibitem[W]{W} P. W. Wong, ``A public key watermark for image verification and authentication," in \textit{Proceedings of IEEE International Conference on Image Processing}, Chicago, USA, October 4-7, 1998, 425-429.
\bibitem[Za]{Za} P. L. Zador, ``Asymptotic quantization error of continuous signals and their quantization dimension," IEEE Trans. Inform. Theory \textbf{28}, 139-148 (1982).
\bibitem[Zw]{Zw} M. Zwick, B. Lovell and J. Marsh, ``Global optimization studies on the 1-D phase problem," International Journal of General Systems \textbf{25}, 47-59 (1996).
\end{thebibliography}
\end{document}